\documentclass[11pt,fleqn,twoside]{article}
\usepackage{amsfonts,amssymb,latexsym}
\makeatletter
\newcommand{\prava}[1]{\small\it
\begin{flushleft}
Copyright \copyright \ 1999 by  #1
\end{flushleft}}

\newcommand{\name}[1]{\begin{flushleft}
                       \LARGE \bf #1
                       \end{flushleft}\vspace{-3mm}}

\newcommand{\Author}[1]{\begin{flushleft}
                       \it #1 \end{flushleft}}

\newcommand{\Adress}[1]{\begin{flushleft}
                       \it #1 \end{flushleft}}

\newcommand{\Date}[1]{\begin{flushleft}
                      \small  \it #1 \end{flushleft}}

\newcommand{\ehkol}{Author \ name}
\newcommand{\ohkol}{Article \ name}
\renewcommand{\@evenhead}{
\hspace*{-3pt}\raisebox{-15pt}[\headheight][0pt]{\vbox{\hbox to \textwidth 
{\thepage \hfil \ehkol}\vskip4pt \hrule}}}
\renewcommand{\@oddhead}{
\hspace*{-3pt}\raisebox{-15pt}[\headheight][0pt]{\vbox{\hbox to \textwidth 
{\ohkol \hfil \thepage}\vskip4pt\hrule}}}
\renewcommand{\@evenfoot}{}
\renewcommand{\@oddfoot}{}

\newcommand{\be}{\begin{equation}}
\newcommand{\ee}{\end{equation}}
\newcommand{\ba}{\hspace*{-5pt}\begin{array}}
\newcommand{\ea}{\end{array}}

\newcommand{\ds}{\displaystyle}
\makeatother

\def\bimp{\mbox{\mathversion {bold}$ P$}}
\def\bimpi{\mbox{\mathversion {bold}$ P$}^{-1}}
\def\bimag{\mbox{\mathversion {bold}$ A$}_{\gamma}}
\def\bima{\mbox{\mathversion {bold}$ A$}}
\def\bimdg{\mbox{\mathversion {bold}$ D$}_{\gamma}}
\def\vfgz{\vec{f}_{\gamma}(z)}
\def\vfg{\vec{f}_{\gamma}}
\def\nvfgz{\|\vfgz\|}

\begin{document}

\thispagestyle{empty}
\setcounter{page}{139}

\renewcommand{\ehkol}{S. Melkonian}
\renewcommand{\ohkol}{Psi-Series Solutions of the Cubic H\'{e}non-Heiles System and
Their Convergence}

\begin{flushleft}
\footnotesize \sf
Journal of Nonlinear Mathematical Physics \qquad 1999, V.6, N~2,
\pageref{melkonian-fp}--\pageref{melkonian-lp}.
\hfill {\sc Article}
\end{flushleft}

\vspace{-5mm}

\renewcommand{\footnoterule}{}
{\renewcommand{\thefootnote}{} 
 \footnote{\prava{S. Melkonian}}}

\name{Psi-Series Solutions of the Cubic H\'{e}non-Heiles System and
Their Convergence}\label{melkonian-fp}

\Author{S. MELKONIAN}

\Adress{School of Mathematics and Statistics, Carleton University,\\
Ottawa, Ontario K1S 5B6, Canada\\[1mm]
E-mail: melkonia@math.carleton.ca}

\Date{Received August 31, 1998; Revised October 22, 1998,
Accepted December 12, 1998}

\begin{abstract}
\noindent
The cubic H\'{e}non-Heiles system contains parameters,
for most values of which, the system is not integrable. In such
parameter regimes, the general solution is expressible in formal
expansions about arbitrary movable branch points, the so-called
psi-series expansions. In this paper, the convergence of known, as
well as new, psi-series solutions on real time intervals is proved,
thereby establishing that the formal solutions are actual
solutions.
\end{abstract}

\renewcommand{\theequation}{\arabic{section}.\arabic{equation}}

\section*{1. Introduction}

\setcounter{section}{1}

\setcounter{equation}{0}

\noindent
The generalized H\'{e}non-Heiles Hamiltonian is
\be
H=\frac{1}{2}\left(\dot x^2+\dot y^2+Ax^2+By^2\right)+
Dx^2y^{m-2}-\frac{C}{m}y^m,
\ee
where $\dot x=\frac{dx}{dt}$, $\dot y=\frac{dy}{dt}$
and $A$, $B$, $C$
and $D$ are nonzero parameters. It corresponds to the original form
given by H\'{e}non and Heiles [1] in their studies of astronomical
problems when $A=B=C=D=1$ and $m=3$. Here, we consider the cubic form
($m=3$) of (1.1), with the corresponding equations of motion
\be
\ba{l}
\ddot x+Ax+2Dxy=0,
\vspace{1mm}\\
\ddot y+By+Dx^2-Cy^2=0.
\ea
\ee
Several studies have been made regarding the integrability properties
and singularity structures of various versions of (1.2); for example,
Chang, Tabor and Weiss~[2] have analy\-sed~(1.2) with $A=B=1$, whereas
Chang, Greene, Tabor and Weiss [3] have considered the same system
for general~$m$. In keeping with their notation, we rescale the
variables,
\be
x\rightarrow \frac{x}{C},\qquad  y\rightarrow \frac{y}{C},\qquad
\lambda=\frac{D}{C},
\ee
to express (1.2) as
\be
\ba{l}
\ddot x+Ax+2\lambda xy=0,
\vspace{1mm}\\
\ddot y+By+\lambda x^2-y^2=0.
\ea
\ee
The leading-order behaviours of solutions about an arbitrary
singularity, as well as their resonance structures (in the sense of
the ARS algorithm, Ablowitz, Ramani and Segur [4, 5]), have been
derived (see [2, 3]; Chang, Tabor, Weiss and Corliss [6]; Bountis,
Segur and Vivaldi [7]; Weiss [8]; Tabor [9, p. 337]).

We summarize the results brief\/ly, in order to f\/ix notation and to
single out the present objectives. Setting $\tau=t-t_0$, where $t_0$
is an arbitrary constant (the location of the singularity), and
\be
x\sim a\tau^{\alpha},\qquad y\sim b\tau^{\beta},\qquad \tau\rightarrow 0,
\ee
two types of leading-order singular behaviours are found:
\be
\ba{ll}
\mbox{(i) }& \alpha=\beta=-2,
\vspace{1mm}\\
\mbox{(ii) }&  \beta=-2,  \quad Re(\alpha)>\beta.
\ea
\ee
Carrying the analysis further by writing
\be
x\sim a\tau^{\alpha}+p\tau^{\alpha+k_r},\qquad
 y\sim b\tau^{\alpha}+q\tau^{\beta+k_r},
\ee
and determining those values of $k_r$ for which $p$ or $q$ is
undetermined, one f\/inds the corresponding resonances. It turns out
that
\be
\ba{ll}
\mbox{(i)}& \alpha=\beta=-2,\quad k_r=-1, 6, r, \bar{r},
\vspace{2mm}\\
&\ds \mbox{where}\quad r,\bar{r}=\frac{5}{2}\pm\frac{1}{2}
\sqrt{1-24(\frac{1}{\lambda}+1)};
\vspace{3mm}\\
\mbox{(ii)} & \ds Re(\alpha)>-2,\quad  \beta=-2,\quad \alpha,
\bar{\alpha}=\frac{1}{2}\pm\frac{1}{2}\sqrt{1-48\lambda},\quad k_r=-1, 6,  0, r,
\bar{r},
\vspace{3mm}\\
&\mbox{where}\quad  r, \bar{r}=\mp\sqrt{1-48\lambda}.
\ea
\ee
In case (ii), $\alpha$ and $\bar{\alpha}$ are two possible leading orders,
with corresponding resonances $r$ and~$\bar{r}$. This case occurs
only for $\lambda>-\frac{1}{2}$, since $Re(\alpha)>-2$.

The nature of the resonances for each value of $\lambda$ is
summarized as follows:

\medskip

\noindent
{\bf Case(i)}
\[
\ba{lcl}
\ds  \lambda<-\frac{24}{23} & : & r \quad \mbox{and} \quad
\bar{r} \quad \mbox{are complex}\vspace{3mm} \\
 \ds \lambda=-\frac{24}{23} & : & \ds r=\bar{r}=\frac{5}{2}\vspace{3mm}\\
\ds -\frac{24}{23}<\lambda<-\frac{1}{2} &:& r>0,\quad  \bar{r}>0,\quad
r\neq\bar{r}\vspace{3mm}\\
\ds \lambda=-\frac{1}{2}&:& r=5,\quad \bar{r}=0\vspace{3mm}\\
\ds -\frac{1}{2}<\lambda<0 &:& r>0,\quad  \bar{r}<0\vspace{3mm}\\
\lambda>0 &:& r \quad \mbox{and} \quad \bar{r} \quad \mbox{are complex}
\ea
\]
{\bf Case (ii)}
\[
\ba{lcl}
\ds  -\frac{1}{2}<\lambda<0& :& r<0,\quad  \bar{r}>0\vspace{3mm}\\
\ds 0<\lambda<\frac{1}{48} &:& r<0,\quad  \bar{r}>0\vspace{3mm}\\
\ds \lambda=\frac{1}{48} &:& r= \bar{r}=0\vspace{3mm}\\
\ds \lambda>\frac{1}{48}&:& r \quad \mbox{and} \quad \bar{r} \quad
\mbox{are pure imaginary}
\ea
\]
Those cases where a negative resonance is present (in addition to $-1$)
likely correspond to singular, rather than general, solutions (see,
e.g., [9, p.~339]). The signif\/icance of repeated resonances (case~(i),
$\lambda=-\frac{24}{23}$, and case (ii), $\lambda=\frac{1}{48}$) is,
at present, unknown. The case $\lambda=-\frac{1}{2}$ is somewhat
anomalous [9, p.~340] and we shall comment upon it later. In those
cases where there exist four distinct resonances, the usual $-1$ and
three others with nonnegative real parts, general solutions may be
written down in series expansions about the arbitrary singularity
$t_0$. These may be Laurent series, if the system passes the
Painlev\'{e} test~[5] or, otherwise, psi-series (Hille~[10, p.~249]).

It has been found [2] that the system is integrable when $\lambda=-1$
and $A=B$, when $\lambda=-\frac{1}{6}$, and when $\lambda=-\frac{1}{16}$ and
$B=16A$. In the latter case, although the system is integrable, it
passes only the ``weak'' Painlev\'{e} test in that it admits a
rational movable branch point. Outside of these three cases, the
system possesses solutions with movable branch points (logarithmic,
irrational or complex), and a psi-series expansion is required to
represent the general solution about an arbitrary branch point. Such
expansions constitute formal general solutions containing four
arbitrary constants, the remaining coef\/f\/icients being well def\/ined by
a self-consistent recursion relation.

In this paper, we prove the absolute convergence of such series on
real intervals of the form $0<\tau<R$, $R>0$, thus establishing that
the formal solutions are actual solutions. Similar results have been
obtained by Melkonian and Zypchen~[11] for the Lorenz system (see,
e.g.,~[9, p.~344] and Sparrow~[12]), where only logarithmic
psi-series occur. A study of the convergence problem for a class of
second- and third-order ordinary dif\/ferential equations has been made
by Hemmi and Melkonian~[13]. Some of the key ideas within the proofs
may be extracted from a paper by Hille~[14] regarding second-order
quadratic systems. Studies of general systems of partial dif\/ferential
equations which admit WTC (Weiss, Tabor and Carnevale~[15])
expansions containing logarithms have been made by Kichenassamy and
Srinivan [16] and Kichenassamy and Littman~[17], who have proved
convergence by entirely dif\/ferent methods. Logarithmic psi-series
occur for the Cubic H\'{e}non Heiles system~(1.4) when $\lambda=-1$ and
$A\neq B$, or when $\lambda=-\frac{1}{16}$ and $B\neq 16A$, and these may
be dealt with by the methods discussed in [11, 13, 16, 17].

Here, we restrict our attention to the non-logarithmic cases.
Sections~2 and~3 deal with the case(i)-leading orders, involving a
series with complex exponents (Section~2) and one with irrational
exponents (Section 3). The necessary Lemmas are relegated to Appendix~A.
The series which we employ in Sections~2 and~3 are not the same as
the ones given, e.g., in~[3], and we conf\/irm the validity
(self-consistency) of our version in Appendix~B. Section~4 concerns
the case(ii)-leading orders where, once again, we employ a new
series, the validity of which is conf\/irmed in Appendix~C. Concluding
remarks are made in Section~5, including a discussion of the case
$\lambda=-\frac{1}{2}$.

\section*{2. Case(i)-leading orders:
{\mathversion{bold}
$\alpha=\beta=-2$,\\[2mm]
\hspace*{20pt}$k_r=-1, 6,
\frac{5}{2}\pm\frac{1}{2}\sqrt{1-24\left(\frac{1}{\lambda}+1\right)}$}
complex}

\setcounter{section}{2}\setcounter{equation}{0}

\noindent
Here, the resonances $r,
\bar{r}=\frac{5}{2}\pm\frac{1}{2}\sqrt{1-24(\frac{1}{\lambda}+1)}$ are complex,
with $\lambda<-\frac{24}{23}$ or $\lambda>0$.
The general solution may be taken in the form
\be
\ba{l}
\ds x=\sum_{k=0}^{\infty}\sum_{l=0}^{\infty}\sum_{m=0}^{\infty}
a_{klm}\tau^{k-2+rl+\bar{r}m},\vspace{3mm}\\
\ds y=\sum_{k=0}^{\infty}\sum_{l=0}^{\infty}\sum_{m=0}^{\infty}b_{klm}
\tau^{k-2+rl+\bar{r}m}.
\ea
\ee
The self-consistency of (2.1), as well as the presence of four
arbitrary constants, is shown in Appendix~B. It is also shown therein
that (2.1) is equivalent to
\be
\ba{l}
\ds x=\sum_{k=0}^{\infty}\sum_{l=0}^{\infty}
\left[a_{kl}\tau^{k-2+rl}+\bar{a}_{kl}\tau^{k-2+\bar{r}l}\right],
\vspace{3mm}\\
\ds y=\sum_{k=0}^{\infty}\sum_{l=0}^{\infty}
\left[b_{kl}\tau^{k-2+rl}+\bar{b}_{kl}\tau^{k-2+\bar{r}l}\right],
\ea
\ee
and that (2.2) is equivalent to the solution given in [3]. Either
(2.1) or (2.2) may be employed to prove convergence. However, the
advantage of starting with (2.1) is that it is much easier to deal
with a single, triply-indexed series such as (2.1) than a sum of two,
doubly-indexed ones such as (2.2), vis-\`{a}-vis the conf\/irmation of
its validity.

The convergence proof consists of resumming (2.2) into more tractable
forms, and showing that the latter are majorized by series which
converge on an interval, by the ratio test. Thus, let
$\varepsilon=\frac{1}{2}\sqrt{24(\frac{1}{\lambda}+1)-1}>0$, note that $r-1,
\bar{r} -1=\frac{3}{2}\pm i\varepsilon$, and write
\be
\ba{l}
\ds x=\sum_{k=0}^{\infty}\sum_{l=0}^{\infty}
\left[a_{kl}\tau^{k-2+l+\left(\frac{3}{2}+i\varepsilon\right)l}
+\bar{a}_{kl}\tau^{k-2+l+\left(\frac{3}{2}-i\varepsilon\right)l}\right]
\vspace{3mm}\\
\ds \qquad
=\sum_{k=0}^{\infty}\sum_{l=0}^{\infty}\left[(a_{kl}+
\bar{a}_{kl})\cos(\varepsilon l z)+i(a_{kl}-\bar{a}_{kl})\sin(\varepsilon
lz)\right]e^{\frac{3}{2}lz}\tau^{k-2+l},
\ea
\ee
where $z=\ln(\tau)$, $\tau>0$. It follows that
\be
x=\sum_{\gamma=0}^{\infty}f_{\gamma}(z)\tau^{\gamma-2},
\ee
where
\be
f_{\gamma}(z)=\sum_{k+l=\gamma}\left[(a_{kl}+\bar{a}_{kl})\cos(\varepsilon
lz)+i(a_{kl}-\bar{a}_{kl})\sin(\varepsilon lz)\right]e^{\frac{3}{2}lz}
\ee
is a polynomial in $e^{\frac{3}{2}z}$ of degree $\gamma$, with
coef\/f\/icients that are bounded functions of $z$. Similarly,
\be
y=\sum_{\gamma=0}^{\infty}h_{\gamma}(z)\tau^{\gamma-2},
\ee
where
\be
h_{\gamma}(z)=\sum_{k+l=\gamma}\left[(b_{kl}+\bar{b}_{kl})\cos(\varepsilon
lz)+i(b_{kl}-\bar{b}_{kl})\sin(\varepsilon lz)\right]e^{\frac{3}{2}lz}
\ee
is a polynomial of the same type as $f_{\gamma}(z)$.

Let $u(t)=\dot{x}(t)\equiv\frac{dx}{dt}$ and
$v(t)=\dot{y}(t)\equiv\frac{dy}{dt}$, then
\be
\ba{l}
\ds u=\dot x=\sum_{\gamma=0}^{\infty}g_{\gamma}(z)
\tau^{\gamma-3}=\sum_{\gamma=0}^{\infty}\left[(\gamma-2)f_{\gamma}
+f_{\gamma}'\right]\tau^{\gamma-3},\vspace{3mm}\\
\ds v=\dot
y=\sum_{\gamma=0}^{\infty}k_{\gamma}(z)\tau^{\gamma-3}=
\sum_{\gamma=0}^{\infty}\left[(\gamma-2)h_{\gamma}+h_{\gamma}'\right]\tau^{\gamma-3},
\ea
\ee
where $f'_{\gamma}\equiv\frac{d f_{\gamma}}{dz}, h'_{\gamma}\equiv\frac{d
h_{\gamma}}{dz}$, and $g_{\gamma}$ and $k_{\gamma}$ are also polynomials of the
same type as $f_{\gamma}$ and~$h_{\gamma}$.

Express (1.4) as a system of four f\/irst-order equations in $x, u, y,
v$ to obtain
\be
\ba{l}
\dot x=u,
\vspace{1mm}\\
\dot u=-Ax-2\lambda xy,
\vspace{1mm}\\
\dot y=v,
\vspace{1mm}\\
\dot v=-By-\lambda x^2+y^2.
\ea
\ee
Substitution of (2.4), (2.6) and (2.8) into (2.9) gives, for
$\gamma\geq 0$,
\be
\ba{l}
\ds f_{\gamma}'+(\gamma-2)f_{\gamma}-g_{\gamma}=0,
\vspace{3mm}\\
\ds g_{\gamma}'+(\gamma-3)g_{\gamma}+
Af_{\gamma-2}+2\lambda\sum_{\mu=0}^{\gamma}f_{\gamma-\mu}h_{\mu}=0,
\vspace{3mm}\\
h_{\gamma}'+(\gamma-2)h_{\gamma}-k_{\gamma}=0,
\vspace{3mm}\\
\ds k_{\gamma}'+(\gamma-3)k_{\gamma}+Bh_{\gamma-2}+\lambda\sum_{\mu=0}^{\gamma}f_{\gamma-\mu}f_{\mu}
-\sum_{\mu=0}^{\gamma}h_{\gamma-\mu}h_{\mu}=0.
\ea
\ee
For $\gamma=0$, we f\/ind
\be
f_0=\pm\frac{3}{\lambda}\sqrt{\frac{1}{\lambda}+2},\qquad g_0=-2f_0,
\qquad h_0=-\frac{3}{\lambda},\qquad k_0=\frac{6}{\lambda}.
\ee
For $\gamma>0$, (2.10) is expressible as
\be
\vec{f}_{\gamma}'+\bimag\vec{f}_{\gamma}=\vec{F}_{\gamma},
\ee
where
\be
\ba{l}
\ds\vec{f}_{\gamma}=\left(\begin{array}{c}f_{\gamma}\\
g_{\gamma}\\h_{\gamma}\\k_{\gamma}\ea\right),
\qquad
\bimag=\left(\begin{array}{cccc}\gamma-2 & -1 &0&0\vspace{2mm}\\
-6&\gamma-3&2\lambda f_0&0\vspace{2mm}\\
0&0&\gamma-2&-1\vspace{2mm}\\
2\lambda f_0&0&\ds \frac{6}{\lambda}&\gamma-3\ea\right),
\vspace{4mm}\\
\vec{F}_{\gamma}=\left(\begin{array}{c}
0\vspace{3mm}\\
\ds -Af_{\gamma-2}-2\lambda\sum\limits_{\mu=1}^{\gamma-1}f_{\gamma-\mu}h_{\mu}
\vspace{3mm}\\
0\vspace{3mm}\\
\ds -Bh_{\gamma-2}-\lambda\sum\limits_{\mu=1}^{\gamma-1}f_{\gamma-\mu}f_{\mu}
+\sum\limits_{\mu=1}^{\gamma-1}h_{\gamma-\mu}h_{\mu}\ea\right)
\equiv\left(\begin{array}{c}F_{\gamma}\\
G_{\gamma}\\H_{\gamma}\\K_{\gamma}\ea\right).
\ea
\ee
The eigenvalues of $\bimag$ are $\gamma+1, \gamma-6$ and
$\gamma-\left(\frac{5}{2}\pm i\varepsilon\right)$, precisely $\gamma-k_r$ where $k_r$
is a resonance. The matrices $\bimag$ are simultaneously
diagonalizable by a matrix $\bimp$ independent of $\gamma$. $\bimp$ is
the matrix of eigenvectors,
\be
\bimp=\left(\begin{array}{cccc}
-\lambda f_0&-\lambda f_0&\lambda f_0&\lambda f_0
\vspace{3mm}\\
3\lambda f_0&-4\lambda f_0&\ds \lambda f_0\left(\frac{1}{2}+i\varepsilon\right)& \ds \lambda
f_0\left(\frac{1}{2}-i\varepsilon\right)
\vspace{3mm}\\
3&3&\ds 3\left(2+\frac{1}{\lambda}\right)&\ds 3\left(2+\frac{1}{\lambda}\right)
\vspace{3mm}\\
-9&12& \ds 3\left(2+\frac{1}{\lambda}\right)\left(\frac{1}{2}
+i\varepsilon\right)&\ds 3\left(2+\frac{1}{\lambda}\right)\left(\frac{1}{2}-
i\varepsilon\right)\ea\right),
\ee
and $\bimpi\bimag\bimp=\bimdg$, where
\be
\bimdg=\left(\begin{array}{cccc}\gamma+1&0&0&0\vspace{2mm}\\
0&\gamma-6&0&0\vspace{2mm}\\
 0&0&\ds \gamma-\left(\frac{5}{2}+i\varepsilon\right)&0\vspace{2mm}\\
0&0&0&\ds \gamma-\left(\frac{5}{2}-i\varepsilon\right)
\ea\right).
\ee
For $\gamma>6$ (in which case all eigenvalues of $\bimag$ have positive
real parts), the solution of the matrix dif\/ferential equation (2.12)
may be expressed as
\be
\vec{f}_{\gamma}(z)=\int_{-\infty}^z e^{\bimag (x-z)}\vec{F}_{\gamma}(x)dx
=\int_{-\infty}^z\bimp e^{\bimdg (x-z)}\bimpi\vec{F}_{\gamma}(x)dx,
\ee
where the exponential within the integrals denotes the exponential of
a matrix.

Given an $m\times n$ matrix $\bima=(a_{ij})$, def\/ine its norm as
\be
\|\bima\|=\max_{1\leq i\leq m}\left(\sum_{j=1}^n|a_{ij}|\right).
\ee
It follows from (2.16) that
\be
\|\vec{f}_{\gamma}(z)\|\leq \|\bimp\|\|\bimpi\|\int_{-\infty}^{z}e^{(\gamma-6)(x-z)}\|
\vec{F}_{\gamma}(x)\|dx,
\ee
since $\|e^{\bimdg(x-z)}\|=e^{(\gamma-6)(x-z)}$ for $x\leq z$. $\bimp$
depends upon $\lambda$, so that $\|\bimp\|\|\bimpi\|\leq M(\lambda)$, a
constant (independent of $\gamma$). Thus,
\be
\|\vec{f}_{\gamma}(z)\|\leq M(\lambda)
\int_{-\infty}^{z}e^{(\gamma-6)(x-z)}\|\vec{F}_{\gamma}(x)\|dx.
\ee

\noindent
{\bf Theorem 2.1.} {\it There exists $K>0$ such that for all $z<0$ and
$\gamma\geq 1$,}
\be
\|\vec{f}_{\gamma}(z)\|\leq\frac{(2K)^{\gamma}}{\sqrt{\gamma+1}}.
\ee

\noindent
{\bf Proof.} Since $f_{\gamma}(z)$ is a polynomial in
$X=e^{\frac{3}{2}z}$ of degree $\gamma$, with coef\/f\/icients that are
bounded functions of $z$ (linear combinations of $\cos(\varepsilon lz)$ and
$\sin(\varepsilon lz)$),
\be
|f_{\gamma}(z) | \leq
\sum_{m=0}^{\gamma}a_{m}^{(\gamma)}e^{\frac{3}{2}mz}\equiv P_{\gamma}(X),
\ee
with $a_m^{(\gamma)}\geq 0$ for $0\leq m\leq \gamma$. By Lemma A1 of
Appendix A (with $q=1,\ n_{\gamma}=\gamma,\ M=1$ and $p=1$), given $N>1$,
there exists $K_1>0$ such that
\be
|f_{\gamma}(z)|\leq P_{\gamma}(X)\leq
\frac{(K_1+KX)^{\gamma}}{\sqrt{\gamma+1}}
\qquad\mbox{for}\quad  1\leq\gamma\leq N-1.
\ee
Similarly, there exist positive constants $K_2$, $K_3$ and $K_4$
corresponding to the polynomials $g_{\gamma}$, $h_{\gamma}$ and
$k_{\gamma}$, repectively. With $K=\max\limits_{1\leq i\leq 4}\{K_i\}$,
we obtain
\be
\|\vec{f}_{\gamma}(z)\|\leq \frac{(K+Ke^{\frac{3}{2}z})^{\gamma}}
{\sqrt{\gamma+1}}\leq\frac{(2K)^{\gamma}}{\sqrt{\gamma+1}}
\qquad\mbox{for}\quad 1\leq\gamma\leq N-1,\ z<0,
\ee
since $0<e^{\frac{3}{2}z}<1$. The inequality (2.23) constitutes our
inductive hypothesis, and we shall show that (2.20) holds for
$\gamma=N$, thus establishing that it holds for all $N>1$.

By (2.13) and the inductive hypothesis (2.23), we f\/ind that
\be
\ba{l}
\ds |G_N|\leq
|A||f_{N-2}|+2|\lambda|\sum_{\mu=1}^{N-1}|f_{N-\mu}||h_{\mu}|
\vspace{3mm}\\
\ds \phantom{G_N}\leq
|A|\frac{(2K)^{N-2}}{\sqrt{N-1}}+2|\lambda|
\sum_{\mu=1}^{N-1}\frac{(2K)^N}{\sqrt{N-\mu+1}\,\sqrt{\mu+1}}
\vspace{3mm}\\
\ds \phantom{G_N}\leq \left
\{\frac{|A|(2K)^{-2}}{\sqrt{N-1}}+2|\lambda|\pi\right\}(2K)^N
\ea
\ee
by Lemma A2 of Appendix A. Similarly,
\be
|K_N|\leq\left\{\frac{|B|(2K)^{-1}}{\sqrt{N-1}}+(|\lambda|+1)\pi
\right\}(2K)^N.
\ee
Thus,
\be
\|\vec{F}_N\|\leq E(2K)^N,
\ee
where
\be
E\leq\max\{|A|+2|\lambda|\pi, |B|+(|\lambda|+1)\pi\}.
\ee
Thus, from (2.19) and (2.26), we obtain
\be
\ba{l}
\ds \|\vec{f}_{N}(z)\| \leq
 M(\lambda)E\int_{-\infty}^z e^{(N-6)(x-z)}(2K)^Ndx
 \vspace{3mm}\\
\ds \phantom{\|\vec{f}_{N}(z)\|} =\left\{\frac{M(\lambda)E\sqrt{N+1}}{N-6}
\right\}\frac{(2K)^N}{\sqrt{N+1}}\qquad
\mbox{for}\quad N>6.
\ea
\ee
For all suf\/f\/iciently large $N$, $\frac{M(\lambda)E\sqrt{N+1}}{N-6}<1$, so
beginning with such an $N>6$, we obtain
\be
\|\vec{f}_N(z)\|\leq \frac{(2K)^N}{\sqrt{N+1}},
\ee
which completes the proof.

By (2.4),
\be
x=\sum_{\gamma=0}^{\infty}f_{\gamma}(z)\tau^{\gamma-2},
\ee
and by (2.20),
\be
|f_{\gamma}(z)\tau^{\gamma-2}|\leq
\frac{(2K)^{\gamma}}{\sqrt{\gamma+1}}\tau^{\gamma-2},
\ee
and the series
\be
\sum_{\gamma=0}^{\infty}\frac{(2K)^{\gamma}}{\sqrt{\gamma+1}}\tau^{\gamma-2}
\ee
converges absolutely by the ratio test for $\tau<\frac{1}{2K}$. Thus,
the series (2.30) for $x$ converges absolutely for $0<\tau<R$, where
$R$ is at least $\frac{1}{2K}$. Similarly, the series for $\dot x,\ y$
and $\dot y$ converge absolutely for $0<\tau<R$.

\section*{3. Case(i)-leading orders: {\mathversion{bold}$\alpha=\beta=-2$,\\[2mm]
\hspace*{20pt} $k_r=-1, 6, \frac{5}{2}\pm\frac{1}{2}\sqrt{1-24\left(\frac{1}
{\lambda}+1\right)}$} irrational}

\setcounter{section}{3}\setcounter{equation}{0}

\noindent
In the present case, the resonances $r,
\bar{r}=\frac{5}{2}\pm\frac{1}{2}\sqrt{1-24\left(\frac{1}{\lambda}+1\right)}$
are irrational, positive and distinct,
with $-\frac{24}{23}<\lambda<-\frac{1}{2}$. The general solution is given
by (2.2) (or (2.1)) as in the complex-resonance  case of Section 2,
but the details of the resummation are slightly dif\/ferent since,
in the present case, $r<1$ and/or $\bar{r}<1$ is possible.

Let $n$ be the least positive integer such that $\mu_1=r-\frac{1}{n}>0$
and $\mu_2=\bar r -\frac{1}{n}>0$. (Any positive integer satisfying the
latter two conditions will also do.) Then the solution $x$ in (2.2)
may be expressed as
\be
x=\sum_{k=0}^{\infty}\sum_{l=0}^{\infty}
\left[a_{kl}\tau^{k-2+\frac{1}{n}l+\mu_1l}
+\bar{a}_{kl}\tau^{k-2+\frac{1}{n}l+\mu_2l}\right].
\ee
Let $w=\tau^{1/n}$ and $z=\ln(\tau)$, $\tau>0$, then
\be
\ba{l}
\ds x=\sum_{k=0}^{\infty}\sum_{l=0}^{\infty}
\left[a_{kl}w^{nk-2n+l}e^{\mu_1
lz}+\bar{a}_{kl}w^{nk-2n+l}e^{\mu_2 lz}\right]
\vspace{3mm}\\
\ds \qquad =\sum_{\gamma=0}^{\infty}\sum_{nk+l=\gamma}\left(a_{kl}e^{\mu_1
lz}+\bar{a}_{kl}e^{\mu_2 lz}\right)w^{\gamma-2n}
=\sum_{\gamma=0}^{\infty}f_{\gamma}(z)w^{\gamma-2n},
\ea
\ee
where
\be
f_{\gamma}(z)=\sum_{nk+l=\gamma}
\left(a_{kl}e^{\mu_1 lz}+\bar{a}_{kl}e^{\mu_2 lz}\right)
\ee
is a polynomial of degree $\gamma$ in the two variables $e^{\mu_1 z}$
and $e^{\mu_2 z}$. Similarly,
\be
y=\sum_{\gamma=0}^{\infty}h_{\gamma}(z)w^{\gamma-2n},
\ee
where
\be
h_{\gamma}(z)=\sum_{nk+l=\gamma}\left(b_{kl}e^{\mu_1 lz}
+\bar{b}_{kl}e^{\mu_2 lz}\right)
\ee
is a polynomial of the same type as $f_{\gamma}(z)$.

Let
\be
\ba{l}
\ds u=\dot x=\sum_{\gamma=0}^{\infty}g_{\gamma}(z)w^{\gamma-3n}
=\sum_{\gamma=0}^{\infty}\left[\left(\frac{\gamma}{n}-2\right)f_{\gamma}+
f_{\gamma}'\right]w^{\gamma-3n},\vspace{3mm}\\
\ds v=\dot
y=\sum_{\gamma=0}^{\infty}k_{\gamma}(z)w^{\gamma-3n}=\sum_{\gamma=0}^{\infty}
\left[\left(\frac{\gamma}{n}-2\right)h_{\gamma}+h_{\gamma}'\right]w^{\gamma-3n},
\ea
\ee
where, as before, $\dot x=\frac{dx}{dt},\
f_{\gamma}'=\frac{df_{\gamma}}{dz}$, etc., and $g_{\gamma}$ and $k_{\gamma}$
are polynomials of the same type as $f_{\gamma}$ and $h_{\gamma}$.
Substitution of (3.2), (3.4) and (3.6) into the system (2.9) gives,
for $\gamma=0$,
\be
f_0=\pm\frac{3}{\lambda}\sqrt{\frac{1}{\lambda}+2},\qquad
g_0=-2f_0,\qquad  h_0=-\frac{3}{\lambda},\qquad k_0=\frac{6}{\lambda},
\ee
(as (2.11)) and, for $\gamma>0$,
\be
\vec{f}_{\gamma}'+\bimag\vec{f}_{\gamma}=\vec{F}_{\gamma}
\ee
(as (2.12)), where $\vec{F}_{\gamma}$ is again given by (2.13), but now,
\be
\bimag=\left(\begin{array}{cccc}
\ds \frac{\gamma}{n}-2&-1&0&0\vspace{2mm}\\
-6&\ds \frac{\gamma}{n}-3&2\lambda f_0&0\vspace{2mm}\\
0&0&\ds \frac{\gamma}{n}-2&-1\vspace{2mm}\\
2\lambda f_0&0&\ds \frac{6}{\lambda}&\ds \frac{\gamma}{n}-3\ea\right),
\ee
and coincides with the $\bimag$ of (2.13) if $n=1$. The eigenvalues
of $\bimag$ in the present case are
$\frac{\gamma}{n}-k_r$, where $k_r=-1,
6, r, \bar{r}$ are the resonances.

The matrix $\bimp$ which diagonalizes all of the $\bimag$ is the same
as (2.14) which, for the present situation, is expressed as
\be
\bimp=\left(\begin{array}{cccc}
-\lambda f_0 &-\lambda f_0& \lambda f_0& \lambda f_0
\vspace{2mm}\\
3\lambda f_0 &-4\lambda f_0 & \lambda f_0(r-2)&\lambda f_0(\bar{r}-2)
\vspace{2mm}\\
3&3&\ds 3\left(2+\frac{1}{\lambda}\right)&\ds 3\left(2+\frac{1}{\lambda}\right)\vspace{2mm}\\
-9&12&\ds 3\left(2+\frac{1}{\lambda}\right)(r-2)&\ds 3\left(2+\frac{1}{\lambda}\right)
(\bar r-2)\ea\right),
\ee
with $\bimpi\bimag\bimp=\bimdg$ and
\be
\bimdg=\left(\begin{array}{cccc}
\ds \frac{\gamma}{n}+1&0&0&0\vspace{2mm}\\
0&\ds \frac{\gamma}{n}-6&0&0
\vspace{2mm}\\
0&0&\ds \frac{\gamma}{n}-r&0\vspace{2mm}\\
 0&0&0&\ds \frac{\gamma}{n}-\bar r \ea\right).
\ee
As before, $\|\bimp\|\,\|\bimpi\|$ is bounded by a constant $M(\lambda)$,
the solution of the system (3.8) is expressible as
\be
\vec{f}_{\gamma}(z)=\int_{-\infty}^z \bimp e^{\bimdg(x-z)}\bimpi \vec{F}_{\gamma}(x)dx
\ee
for $\frac{\gamma}{n}>6$ (noting that within the range of $\lambda$ presently
under consideration, $0<r$,  $\bar r<6$, so that all eigenvalues of
$\bimag$ are positive for $\frac{\gamma}{n}>6$), and
\be
\|\vec{f}_{\gamma}(z)\|\leq M(\lambda)\int_{-\infty}^z e^{\left(\frac{\gamma}{n}-6\right)
(x-z)}\|\vec{F}_{\gamma}(x)\|dx.
\ee

\noindent
{\bf Theorem 3.1.} {\it There exists $K>0$ such that for all $z<0$
and $\gamma\geq 1$,}
\be
\nvfgz\leq\frac{(3K)^{\gamma}}{\sqrt{\gamma+1}}.
\ee

\noindent
{\bf Proof.} Since $f_{\gamma}(z)$ is a polynomial of degree $\gamma$ in
$X=e^{\mu_1z}$ and $Y=e^{\mu_2z}$, Lemma A3 of Appendix A (with $M=1$
and $p=1$) implies that given $N>1$, there exists $K_1>0$ such that
\be
|f_{\gamma}(z)|\leq
\frac{(K_1+K_1e^{\mu_1z}+K_1e^{\mu_2z})^{\gamma}}{\sqrt{\gamma+1}}
\qquad\mbox{for}\quad 1\leq\gamma\leq N-1.
\ee
Similarly, there exist positive constants $K_2$, $K_3$ and $K_4$
corresponding to the polynomials $g_{\gamma}$, $h_{\gamma}$ and $k_{\gamma}$,
respectively. With $K=\max\limits_{1\leq i\leq 4}\{K_i\}$, we obtain
\be
\nvfgz\leq \frac{(K+Ke^{\mu_1z}+Ke^{\mu_2z})^{\gamma}}{\sqrt{\gamma+1}}
\qquad\mbox{for}\quad 1\leq\gamma\leq N-1.
\ee
Since $0<e^{\mu_i z}<1$ for $i=1, 2$ and for all $z<0$, we have
\be
\nvfgz\leq\frac{(3K)^{\gamma}}{\sqrt{\gamma+1}}\qquad\mbox{for}\quad
1\leq\gamma\leq N-1.
\ee
The inequality (3.17) constitutes our inductive hypothesis, and we
shall show that (3.14) holds for $\gamma=N$, thereby establishing that it
holds for all $N>1$.

As in Section 2, we f\/ind that
\be
\|\vec{F}_N\|\leq E(3K)^N,
\ee
where $E$ is bounded as in (2.27). Thus, from (3.13) and (3.18), we
obtain
\be
\ba{l}
\ds \|\vec{f}_N(z)\|\leq
M(\lambda)E\int_{-\infty}^ze^{\left(\frac{N}{n}-6\right)(x-z)}(3K)^Ndx
\vspace{3mm}\\
\ds \phantom{\|\vec{f}_N(z)\|}=\left\{\frac{M(\lambda)E\sqrt{N+1}}
{\frac{N}{n}-6}\right\}\frac{(3K)^N}{\sqrt{N+1}}
\qquad\mbox{for}\quad \frac{N}{n}>6.
\ea
\ee
For all suf\/f\/iciently large $N$ (and any $n\geq 1$),
$\frac{M(\lambda)E\sqrt{N+1}}{\frac{N}{n}-6}<1$, so beginning with such an
$N>6n$, we obtain
\be
\|\vec{f}_N(z)\|\leq \frac{(3K)^N}{\sqrt{N+1}},
\ee
which completes the proof.

By (3.2),
\be
x=\sum_{\gamma=0}^{\infty}f_{\gamma}(z)w^{\gamma-2n},
\ee
and by (3.14),
\be
|f_{\gamma}(z)w^{\gamma-2n}|\leq\frac{(3K)^{\gamma}}{\sqrt{\gamma+1}}w^{\gamma-2n},
\ee
and the series
\be
\sum_{\gamma=0}^{\infty}\frac{(3K)^{\gamma}}{\sqrt{\gamma+1}}w^{\gamma-2n}
\ee
converges absolutely by the ratio test for $w<\frac{1}{3K}$, i.e., for
$\tau<\frac{1}{(3K)^n}$. Thus, the series (3.21) for $x$ converges
absolutely for $0<\tau<R$, where $R$ is at least $\frac{1}{(3K)^n}$.
Similarly, the series for $\dot x$, $y$ and $\dot y$ converge
absolutely for $0<\tau<R$.

\section*{4. Case(ii)-leading orders:
{\mathversion{bold} $Re(\alpha)>-2$, $\beta=-2$}}

\setcounter{section}{4}\setcounter{equation}{0}

\noindent
As mentioned in Section 1, if $Re(\alpha)>-2$ and $\lambda>-\frac{1}{2}$, two
types of leading-order and resonance structures are possible, namely,
\be
\ba{l}
\ds \mbox{(a) Leading order}\quad
\alpha=\frac{1}{2}+\frac{1}{2}\sqrt{1-48\lambda}\,,\ \mbox{resonances}\quad
k_r=-1, 6, 0, r;\vspace{3mm}\\
\ds \mbox{(b) Leading order}\quad
\bar\alpha=\frac{1}{2}-\frac{1}{2}\sqrt{1-48\lambda}\,,\ \mbox{resonances}\quad
k_r=-1, 6, 0, \bar r,
\ea
\ee
where $r, \bar r=\mp\sqrt{1-48\lambda}$. For $\lambda<\frac{1}{48}$, $r<0$, so
that (a) does not correspond to a general solution, and will
therefore not be considered. (See Conte, Fordy and Pickering [18] for
a discussion of negative resonances and the appropriate series
expansions.) In case (b), $\bar r>0$ for $-\frac{1}{2}<\lambda<0$ and
$0<\lambda<\frac{1}{48}$, and $\bar r$ is pure imaginary for
$\lambda>\frac{1}{48}$. The general solution takes the form
\be
\ba{l}
\ds x=\sum_{k=0}^{\infty}\sum_{l=0}^{\infty}\sum_{m=0}^{\infty}
a_{klm}\tau^{k+\bar\alpha +(2+\bar\alpha)m+\bar r l},
\vspace{3mm}\\
\ds y=\sum_{k=0}^{\infty}\sum_{l=0}^{\infty}\sum_{m=0}^{\infty}
b_{klm}\tau^{k-2+(2+\bar\alpha)m+\bar r l},
\ea
\ee
discussion of the validity of which, as well as other relevant
issues, is relegated to Appendix~C. If $\bar r$ is pure imaginary
($\lambda>\frac{1}{48}$), then the methods employed earlier do not apply to
(4.2) (see Section~5), and we have not determined whether the series
converge or not. For $\bar r>0$ ($\lambda<\frac{1}{48}$), we proceed as in
the preceding sections.

Let $n$ be the least positive integer such that $\bar\beta=\bar
r-\frac{1}{n}>0$, let $w=\tau^{1/n}$ and $z=\ln(\tau),\ \tau>0$, in
order to express $x$ in (4.2) as
\be
\ba{l}
\ds x=\sum_{k=0}^{\infty}\sum_{l=0}^{\infty}
\sum_{m=0}^{\infty}a_{klm}\tau^{k+2m+\frac{1}{n}l+\bar\alpha+\bar\alpha
m+\bar\beta l}
\vspace{3mm}\\
\ds \qquad =\sum_{k=0}^{\infty}\sum_{l=0}^{\infty}\sum_{m=0}^{\infty}
a_{klm}w^{nk+2nm+l+\bar\alpha n}\tau^{\bar\alpha m+\bar\beta l}
=\sum_{\gamma=0}^{\infty}f_{\gamma}(z)w^{\gamma+\bar\alpha n},
\ea
\ee
where
\be
f_{\gamma}(z)=\sum_{nk+2nm+l=\gamma}a_{klm}e^{\bar\alpha mz+\bar\beta lz}
\ee
is a polynomial of degree $\gamma$ in the two variables $e^{\bar\alpha z}$
and $e^{\bar\beta z}$. Similarly,
\be
y=\sum_{\gamma=0}^{\infty}h_{\gamma}(z)w^{\gamma-2n},
\ee
where
\be
h_{\gamma}(z)=\sum_{nk+2nm+l=\gamma}b_{klm}e^{\bar\alpha mz+\bar\beta lz}
\ee
is a polynomial of the same type as $f_{\gamma}(z)$.

Let
\be
\ba{l}
\ds u=\dot x=\sum_{\gamma=0}^{\infty}g_{\gamma}(z)w^{\gamma+\bar\alpha
n-n}=\sum_{\gamma=0}^{\infty}\left[\left(\frac{\gamma}{n}+\bar
\alpha\right)f_{\gamma}+f_{\gamma}'\right]w^{\gamma+\bar\alpha n-n},
\vspace{3mm}\\
\ds v=\dot y=\sum_{\gamma=0}^{\infty}k_{\gamma}(z)w^{\gamma-3n}
=\sum_{\gamma=0}^{\infty}\left[\left(\frac{\gamma}{n}-2\right)h_{\gamma}+h_{\gamma}'
\right]w^{\gamma-3n},
\ea
\ee
and substitute (4.3), (4.5) and (4.7) into the system (2.9) to
obtain, for $\gamma=0$,
\be
f_0\ \mbox{ arbitrary},\qquad g_0=\bar\alpha f_0,\qquad  h_0=6,
\qquad k_0=-12,
\ee
and for $\gamma>0$, the system
\be
\vfg'+\bimag\vfg=\vec{F}_{\gamma},
\ee
where, as before,
\be
\vfg=\left(\begin{array}{c} f_{\gamma}\\
g_{\gamma}\\h_{\gamma}\\k_{\gamma}\ea\right),\qquad
\vec{F}_{\gamma}=\left(\begin{array}{c} F_{\gamma}\\
G_{\gamma}\\H_{\gamma}\\K_{\gamma}\ea\right),
\qquad F_{\gamma}=H_{\gamma}=0,
\ee
and here,
\be
\ba{l}
\ds G_{\gamma}=-Af_{\gamma-2n}-2\lambda\sum_{\mu=1}^{\gamma-1}f_{\gamma-\mu}h_{\mu},
\vspace{3mm}\\
\ds K_{\gamma}=-Bh_{\gamma-2n}-\lambda\sum_{\mu=0}^{\gamma-4n}e^{2\bar\alpha
z}f_{\gamma-4n-\mu}f_{\mu}+\sum_{\mu=1}^{\gamma-1}h_{\gamma-\mu}h_{\mu},
\ea
\ee
and
\be
\bimag=\left(\begin{array}{cccc}
\ds \frac{\gamma}{n}+\bar\alpha & -1 & 0&0\vspace{2mm}\\
12\lambda & \ds \frac{\gamma}{n}+\bar\alpha-1&2\lambda f_0&0\vspace{2mm}\\
0&0&\ds \frac{\gamma}{n}-2&-1\vspace{2mm}\\
0&0&-12&\ds \frac{\gamma}{n}-3\ea\right).
\ee

As expected, the eigenvalues of $\bimag$ are $\frac{\gamma}{n}-k_r$, where
$k_r=-1, 0, 6, \bar r$ are the resonances. The diagonalizing matrix
for $\bimag$ is
\be
\bimp=\left(\begin{array}{cccc}
\ds -\frac{\bar \alpha}{12}f_0&1&-\lambda f_0&1\vspace{2mm}\\
\lambda f_0&\bar\alpha&-\lambda f_0(\bar\alpha+6)&1-\bar\alpha\vspace{2mm}\\
1&0&3(2\bar\alpha+5)&0\vspace{2mm}\\
-3&0&12(2\bar\alpha+5)&0\ea\right),
\ee
with $\bimdg=\bimpi\bimag\bimp$, and $\|\bimp\|\,\|\bimpi\|\leq
M(\lambda,f_0)$, a constant independent of $\gamma$, but depending upon $\lambda$
and the arbitrary constant $f_0$. Since $-\frac{1}{2}<\lambda<0$ or
$0<\lambda<\frac{1}{48}$, we have $0<\bar r<5<6$, so that $k_r=6$ is again
the largest resonance, and all eigenvalues of $\bimag$ are positive
for $\frac{\gamma}{n}>6$. Thus, for $\frac{\gamma}{n}>6$, the solution of (4.9)
may be expressed as
\be
\vfgz=\int_{-\infty}^z \bimp e^{\bimdg(x-z)}\bimpi\vec{F}_{\gamma}(x)dx,
\ee
with
\be
\nvfgz\leq M(\lambda,f_0)\int_{-\infty}^z e^{\left(\frac{\gamma}{n}-6\right)
(x-z)}\|\vec{F}_{\gamma}(x)\|dx.
\ee

\noindent
{\bf Theorem 4.1.} {\it There exists $K>0$ such that for all $z<0$
and $\gamma\geq 1$},
\be
\nvfgz\leq \frac{(3K)^{\gamma}}{\sqrt{\gamma+1}}.
\ee

\noindent
As the proof is practically identical to the one of Theorem 3.1, we
omit the details.

By (4.3),
\be
x=\sum_{\gamma=0}^{\infty}f_{\gamma}(z)w^{\gamma+\bar\alpha n}
\ee
and by (4.16),
\be
\left|f_{\gamma}(z)w^{\gamma+\bar\alpha n}\right|\leq
\frac{(3K)^{\gamma}}{\sqrt{\gamma+1}}w^{\gamma+\bar\alpha n},
\ee
and the series
\be
\sum_{\gamma=0}^{\infty}\frac{(3K)^{\gamma}}{\sqrt{\gamma+1}}w^{\gamma+\bar\alpha n}
\ee
converges absolutely by the ratio test for $w<\frac{1}{3K}$, i.e., for
$\tau<\frac{1}{(3K)^n}$. Thus, the series (4.17) for $x$ and,
similarly, the series (4.5) for $y$ and (4.7) for $\dot x$ and $\dot
y$ converge absolutely for $0<\tau<R$, where $R$ is at least
$\frac{1}{(3K)^n}$.

\section*{5. Concluding Remarks}

 \setcounter{section}{5}\setcounter{equation}{0}

\noindent
Three types of psi-series solutions of the cubic H\'{e}non-Heiles
system have been shown to be absolutely convergent on intervals of
the form $0<\tau<R$, $R>0$. It is clear that the series for $x$, $\dot x$,
$y$ and $\dot y$ converge uniformly on compact subintervals of
$(0,R)$, thereby justifying the termwise dif\/ferentiation of the
series.

The resummations that have been performed in order to proceed with
the convergence proofs amount to the replacement of multiply-indexed
series with constant coef\/f\/icients by single series with polynomial
coef\/f\/icients. This procedure makes use of the fact that the
resonances appearing as exponents within the psi-series have positive
real parts. For example, to obtain (2.3) and (2.6) from (2.2), it is
essential that $Re(r)>0$ and $Re(\bar{r})>0$. Otherwise, as occurs in
(4.2) when $\bar{r}$ is pure imaginary, a resummation would give rise
to series containing $f_{\gamma}(z)$ which are not polynomials, but
inf\/inite series.

The situation when $\lambda=-\frac{1}{2}$ is very interesting. Here,
$\alpha=-2$ (the leading order for $x$), $r=5$ and $\bar r=0$. But the
leading coef\/f\/icient for $x$ is $\pm\frac{3}{\lambda}\sqrt{2+\frac{1}{\lambda}}=0$,
contradicting the fact that $x=O(\tau^{-2})$. The correct
leading-order behaviours have been given in [7] and~[9, p.~340]. In
our notation,
\[
\ba{l}
\ds x\,\sim\,(-30)^{1/2}\tau^{-2}(\ln(\tau))^{-\frac{1}{2}},
\vspace{3mm}\\
\ds y\,\sim\,6\tau^{-2}+\frac{5}{2}\tau^{-2}(\ln(\tau))^{-1}.
\ea
\]
A full series expansion of the form
\[
\ba{l}
\ds x=\sum_{k=0}^{\infty}\sum_{l=0}^{\infty}a_{kl}\tau^{k-2}
(\ln(\tau))^{-l-\frac{1}{2}},\vspace{3mm}\\
\ds y=\sum_{k=0}^{\infty}\sum_{l=0}^{\infty}b_{kl}\tau^{k-2}(\ln(\tau))^{-l},
\ea
\]
fails, due to an incompatible resonance at $(k,l)=(0,2)$. We are not
aware of the correct series in this case.

\subsection*{Acknowledgments}

\noindent
This work was supported in part by the Natural Sciences and
Engineering Research Council of Canada. I wish to thank the referee
for helpful comments and references.

\renewcommand{\theequation}{\Alph{section}\arabic{equation}}
\section*{Appendix A}

\setcounter{equation}{0}

\noindent
This section contains the Lemmas necessary for the convergence
proofs.

\medskip

\noindent
{\bf Lemma A1.} {\it Let $X\neq 0$ have constant sign, let $q\geq 1$
be an integer, and for $\gamma\geq1$, let
\[
P_{\gamma}(X)=\sum_{m=0}^{n_{\gamma}}c_m^{(\gamma)}X^m
\]
be a sequence of
polynomials of degree $n_{\gamma}=[\frac{\gamma}{q}]$. Given an integer $N>q$,
$0<M\leq 1$ and $p>0$, there exists $K$ with}
$\mbox{sign}(K)=\mbox{sign}(X)$ {\it and} $p|K|>1$ {\it such that}
\[
|P_{\gamma}(X)|\leq M\frac{(p|K|+KX)^{\gamma/q}}{\sqrt{\gamma+1}}\qquad
\mbox{\it for}\quad 1\leq \gamma\leq N-1.
\]
This is Lemma 4.4 of [13], to which we
refer the reader for a proof.

\medskip

\noindent
{\bf Lemma A2.}
\[
\lim_{\gamma \rightarrow\infty}\sum_{\mu=0}^{\gamma-1}\frac{1}{\sqrt{\mu+1}\,
\sqrt{\gamma-\mu}}=\pi.
\]

\noindent
The result follows easily by the proof of the integral test for
convergence of series. For the details, we refer to Hemmi [19, p.~64].

\medskip

\noindent
{\bf Lemma A3.} {\it Let $X, Y>0$ and for $\gamma\geq 1$, let
$P_{\gamma}(X,Y)$ be a sequence of polynomials of degree $\gamma$ in $X$ and
$Y$, i.e.,
\[
P_{\gamma}(X,Y)=\sum_{\beta=0}^{\gamma}\sum_{\mu=0}^{\beta}c_{\mu\beta}^{(\gamma)}X^{\mu}
Y^{\beta-\mu}.
\]
Given an integer $N>1$, $0<M\leq 1$ and $p>0$, there exists $K>0$
such that
\[
|P_{\gamma}(X,Y)|\leq M\frac{(pK+KX+KY)^{\gamma}}{\sqrt{\gamma+1}}
\qquad\mbox{for}\quad  1\leq \gamma\leq N-1.
\] }

\noindent
{\bf Proof.} Let
\[
K=\max_{\mbox{\scriptsize $\begin{array}{c}1\leq\gamma\leq
N-1\\0\leq\beta\leq\gamma\\0\leq\mu\leq\beta\ea$}}
\left\{\frac{|c_{\mu\beta}^{(\gamma)}|
\sqrt{\gamma+1}}{\left(\begin{array}{c}\gamma\\ \beta\ea\right)\left(
\begin{array}{c}\beta\\ \mu\ea\right)Mp^{\gamma-\beta}}\right\}^{1/\gamma},
\]
then for $1\leq\gamma\leq N-1$, $0\leq\beta\leq\gamma$,
$0\leq\mu\leq\beta$, we have
\[
|c_{\mu\beta}^{(\gamma)}|\leq\frac{M}{\sqrt{\gamma+1}}\left(\begin{array}{c}\gamma
\\ \beta\ea\right)\left(\begin{array}{c}\beta \\ \mu\ea\right)K^{\gamma}
p^{\gamma-\beta},
\]
which implies that
\[
\ba{l}
\ds |P_{\gamma}(X,Y)|\leq\sum_{\beta=0}^{\gamma}\sum_{\mu=0}^{\beta}
|c_{\mu\beta}^{(\gamma)}|X^{\mu}Y^{\beta-\mu}
\vspace{3mm}\\
\ds \phantom{|P_{\gamma}(X,Y)|} \leq\frac{M}{\sqrt{\gamma+1}}\sum_{\beta=0}^{\gamma}
\sum_{\mu=0}^{\beta} \left(\begin{array}{c}\gamma \\ \beta\ea\right)
\left(\begin{array}{c}\beta \\ \mu\ea\right)K^{\gamma}p^{\gamma-\beta}X^{\mu}
Y^{\beta-\mu}\vspace{3mm}\\
\ds \phantom{|P_{\gamma}(X,Y)|}=
\frac{M}{\sqrt{\gamma+1}}\sum_{\beta=0}^{\gamma}\sum_{\mu=0}^{\beta}
\left(\begin{array}{c}\gamma \\ \beta\ea\right)\left(\begin{array}{c}
\beta \\ \mu\ea\right) (pK)^{\gamma-\beta}(KX)^{\mu}(KY)^{\beta-\mu}
\vspace{3mm}\\
\ds \phantom{|P_{\gamma}(X,Y)|}=M\frac{(pK+KX+KY)^{\gamma}}{\sqrt{\gamma+1}},
\ea
\]
by the binomial theorem.

\section*{Appendix B}

\setcounter{section}{2}
\setcounter{equation}{0}

\noindent
First, we show that the solution (2.1), namely,
\be
\ba{l}
\ds x=\sum_{k=0}^{\infty}\sum_{l=0}^{\infty}
\sum_{m=0}^{\infty}a_{klm}\tau^{k-2+rl+\bar r m},
\vspace{3mm}\\
\ds y=\sum_{k=0}^{\infty}\sum_{l=0}^{\infty}
\sum_{m=0}^{\infty}b_{klm}\tau^{k-2+rl+\bar r m},
\ea
\ee
of the system (1.4), corresponding to the case(i)-leading orders and
non-integral resonances, is self-consistent and contains four
arbitrary constants.

Substitution of (B1) into (1.4) gives
\be
a_{000}=\pm\frac{3}{\lambda}\sqrt{2+\frac{1}{\lambda}},
\qquad b_{000}=-\frac{3}{\lambda},
\ee
and for $(k,l,m)\neq (0,0,0)$, the coef\/f\/icient recursion relations
\be
\ba{l}
\ds (k+rl+\bar r m)(k+rl+\bar rm -5)a_{klm}+2\lambda a_{000}b_{klm}
\vspace{3mm}\\
\ds \qquad=-Aa_{k-2,l,m}-2\lambda\sum_{(0,0,0)\prec(p,q,s)\prec(k,l,m)}
a_{pqs}b_{k-p,l-q,m-s},
\vspace{3mm}\\
\ds 2\lambda a_{000}a_{klm}+\left[(k+rl+\bar r m -2)(k+rl+\bar r
m-3)+\frac{6}{\lambda}\right]b_{klm}
\vspace{3mm}\\
\ds \qquad=-B b_{k-2,l,m}+
\hspace{-0.7cm}\sum_{(0,0,0)\prec(p,q,s)\prec(k,l,m)}
\hspace{-0.7cm}[-\lambda a_{pqs}a_{k-p,l-q,m-s}+b_{pqs}
b_{k-p,l-q,m-s}],
\ea
\ee
where we have employed the notation $(0,0,0)\prec(p,q,s)\prec(k,l,m)$
to mean that $0\leq p\leq k$, $0\leq q\leq l$, $0\leq s\leq m$, with
$(p,q,s)\neq (0,0,0)$ and $(p,q,s)\neq(k,l,m)$. Expressed in matrix
form, the left-hand sides become
\be\hspace*{-10mm}
\left[\!
\begin{array}{cc}
(k+rl+\bar r m)(k+rl+\bar r m-5) \!\!
& 2\lambda a_{000}\vspace{2mm}\\
 2\lambda a_{000}&  \!\!\! (k+rl+\bar r m-2)(k+rl+\bar r
m-3)+\frac{6}{\lambda}\ea\!\right]\!\!
\left[\begin{array}{c}a_{klm}\\b_{klm}\ea\right],
\ee
and the coef\/f\/icient matrix is singular precisely when $k+rl+\bar r m$
has one of the values $-1, 6, r, \bar r$, i.e., when
\be
(k,l,m)=(-1,0,0), (6,0,0), (0,1,0), (0,0,1).
\ee
These are the resonances of (B1), and compatibility must be checked
at these values, conf\/irming both the self-consistency of (B1) and the
presence of four arbitrary constants.

The resonance at $(k,l,m)=(-1,0,0)$ is trivially compatible, and
corresponds to the arbitrariness of $t_0$. At $(k,l,m)=(0,1,0)$, we
f\/ind that
\be
a_{010}\quad  \mbox{is arbitrary,} \qquad
b_{010}=-\frac{r(r-5)}{2\lambda a_{000}}a_{010},
\ee
and at $(k,l,m)=(0,0,1)$, that
\be
a_{001}\quad \mbox{is arbitrary,} \qquad
b_{001}=-\frac{\bar r(\bar r-5)}{2\lambda
a_{000}}a_{001}.
\ee
At $(k,l,m)=(6,0,0)$, the compatibility condition is
\be
\lambda a_{000}a_{400}A+2\lambda^2a_{000}(a_{200}b_{400}
+a_{400}b_{200})-3b_{400}B  -6(\lambda a_{200}a_{400}-b_{200}b_{400})=0.
\ee
The coef\/f\/icients which af\/fect (B8) are
\be
\ba{l}
\ds a_{000}=\pm\frac{3}{\lambda}\sqrt{2+\frac{1}{\lambda}},\qquad
 b_{000}=-\frac{3}{\lambda},
 \vspace{3mm}\\
 \ds a_{100}=a_{300}=a_{500}=0,\qquad b_{100}=b_{300}=b_{500}=0,
 \vspace{3mm}\\
\ds a_{200}=\frac{a_{000}\left(\frac{A}{\lambda}+B\right)}
{12\left(1+\frac{1}{\lambda}\right)},\qquad
b_{200}=\frac{B-A\left(2+\frac{1}{\lambda}\right)}
{4\lambda\left(1+\frac{1}{\lambda}\right)},
\vspace{3mm}\\
\ds a_{400}=\frac{a_{200}\left(1+\frac{3}{\lambda}\right)
(A+2\lambda b_{200})-\lambda a_{000}\left(Bb_{200}+\lambda a^2_{200}-
b^2_{200}\right)} {10\left(4+\frac{3}{\lambda}\right)},
\vspace{3mm}\\
\ds b_{400}=-\frac{\lambda a_{000}a_{200}(A+2\lambda b_{200})
-2\left(Bb_{200}+\lambda a^2_{200}-b^2_{200}\right)}
{10\left(4+\frac{3}{\lambda}\right)},
\ea
\ee
from which (B8) is conf\/irmed. The system governing $a_{600}$ and
$b_{600}$ is
\be
\left(\begin{array}{cc}6&2\lambda a_{000}
\vspace{2mm}\\
2\lambda a_{000}&
12+\frac{6}{\lambda}\ea\right)
\left(\begin{array}{c}a_{600}\\
b_{600}\ea\right)= \left(\begin{array}{c}-Aa_{400}-2\lambda(a_{200}b_{400}
+a_{400}b_{200})\vspace{3mm}\\
-Bb_{400}-2(\lambda a_{200}a_{400}-b_{200}b_{400})\ea\right),
\ee
and reduces to
\be
6a_{600}+2\lambda a_{000}b_{600}=-A
a_{400}-2\lambda(a_{200}b_{400}+a_{400}b_{200}),
\ee
showing that $a_{600}$ is arbitrary, and $b_{600}$ is def\/ined in
terms of $a_{600}$.

Next, we demonstrate that (2.1) is equivalent to (2.2). The series
which def\/ines $x$ in (2.1) may be broken up into three parts: $m=l$,
$m<l$, and $m>l$. Thus,
\be
\ba{l}
\ds x=\sum_{k=0}^{\infty}\sum_{l=0}^{\infty}a_{kll}\tau^{k-2+5l}
+\sum_{k=0}^{\infty}\sum_{l=1}^{\infty}
\sum_{m=0}^{l-1}a_{klm}\tau^{k-2+5m+(l-m)r}
\vspace{3mm}\\
\ds \qquad+\sum_{k=0}^{\infty}\sum_{m=1}^{\infty}
\sum_{l=0}^{m-1}a_{klm}\tau^{k-2+5l+(m-l)\bar r},
\ea
\ee
where use has been made of the fact that $r+\bar r=5$. The f\/irst sum
in (B12) is expressible as
\be
\sum_{k=0}^{\infty}\sum_{l=0}^{\infty}a_{kll}\tau^{k-2+5l}=
\sum_{i=0}^{\infty}\sum_{k+5l=i}a_{kll}\tau^{i-2}
=\sum_{i=0}^{\infty}a_{i0}\tau^{i-2},
\ee
where
\be
a_{i0}=\sum_{k+5l=i}a_{kll}.
\ee
The second sum in (B12) may be manipulated as follows:
\be
\ba{l}
\ds \sum_{k=0}^{\infty}\sum_{l=1}^{\infty}\sum_{m=0}^{l-1}a_{klm}
\tau^{k-2+5m+(l-m)r}\vspace{3mm}\\
\ds \qquad =\sum_{k=0}^{\infty}\sum_{l=1}^{\infty}
\sum_{j=1}^{l}a_{k,l,l-j}\tau^{k-2+5(l-j)+jr}\qquad (\mbox{letting}\quad
j=l-m)\vspace{3mm}\\
\ds \qquad =\sum_{k=0}^{\infty}\sum_{j=1}^{\infty}\sum_{l=j}^{\infty}
a_{k,l,l-j}\tau^{k-2+5(l-j)+jr}\qquad(\mbox{reversing the order of summation})
\vspace{3mm}\\
\ds \qquad =\sum_{k=0}^{\infty}\sum_{j=1}^{\infty}
\sum_{l=0}^{\infty}a_{k,l+j,l}\tau^{k-2+5l+jr}\qquad(\mbox{replacing $l$
by $l+j$})\vspace{3mm}\\
\ds\qquad  =\sum_{i=0}^{\infty}\sum_{j=1}^{\infty}
\sum_{k+5l=i}a_{k,l+j,l}\tau^{i-2+jr}\vspace{3mm}\\
\ds \qquad  =\sum_{i=0}^{\infty}\sum_{j=1}^{\infty}a_{ij}\tau^{i-2+jr},
\ea\hspace{-7.6pt}
\ee
where
\be
a_{ij}=\sum_{k+5l=i}a_{k,l+j,l}.
\ee
Similarly, the third term in (B12) becomes
\be
\sum_{i=0}^{\infty}\sum_{j=1}^{\infty}\bar a_{ij}\tau^{i-2+j\bar r},
\ee
where
\be
\bar a_{ij}=\sum_{k+5m=i}a_{k,m,m+j}.
\ee
Combining (B13), (B15) and (B17), we obtain
\be
x=\sum_{i=0}^{\infty}\sum_{j=0}^{\infty}a_{ij}\tau^{i-2+jr}
+\sum_{i=0}^{\infty}\sum_{j=1}^{\infty}\bar
a_{ij}\tau^{i-2+j\bar r}.
\ee
Similarly,
\be
y=\sum_{i=0}^{\infty}\sum_{j=0}^{\infty}b_{ij}\tau^{i-2+jr}
+\sum_{i=0}^{\infty}\sum_{j=1}^{\infty}\bar
b_{ij}\tau^{i-2+j\bar r}.
\ee
Def\/ining $\bar a_{i0}=\bar b_{i0}=0$ for all $i\geq 0$, (B19), (B20)
coincide with (2.2).

Finally, the solution given in [3] is
\be
x=\sum_{j=0}^{\infty}\sum_{i=0}^{\infty}A_{ji}\tau^{i-2+j(r-2)}+
\sum_{j=1}^{\infty}\sum_{i=0}^{\infty}\bar
A_{ji}\tau^{i-2+j(\bar r-2)}
\ee
(their equation (3.5a), in our notation), with a similar expression
for $y$. To obtain the form (B21) from (B19), write
\be
\ba{l}
\ds
x=\sum_{i=0}^{\infty}\sum_{j=0}^{\infty}a_{ij}\tau^{i-2+jr}
+\sum_{i=0}^{\infty}\sum_{j=1}^{\infty}\bar
a_{ij}\tau^{i-2+j\bar r}
\vspace{3mm}\\
\ds \qquad =\sum_{i=0}^{\infty}\sum_{j=0}^{\infty}a_{ij}\tau^{i+2j-2+j(r-2)}
+\sum_{i=0}^{\infty}\sum_{j=1}^{\infty}\bar
a_{ij}\tau^{i+2j-2+j(\bar r-2)}
\vspace{3mm}\\
\ds \qquad =\sum_{j=0}^{\infty}\sum_{i=2j}^{\infty}a_{i-2j,j}\tau^{i-2+j(r-2)}
+\sum_{j=1}^{\infty}\sum_{i=2j}^{\infty}\bar
a_{i-2j,j}\tau^{i-2+j(\bar r-2)}
\vspace{3mm}\\
\ds \qquad =\sum_{j=0}^{\infty}\sum_{i=0}^{\infty}A_{ji}\tau^{i-2+j(r-2)}
+\sum_{j=1}^{\infty}\sum_{i=0}^{\infty}\bar
A_{ji}\tau^{i-2+j(\bar r-2)},
\ea
\ee
where
\be
A_{ji}=\left\{\begin{array}{ll}a_{i-2j,j},&\mbox{if }\ i\geq
2j\\0,&\mbox{otherwise}\ea\right\} ,\qquad  \bar
A_{ji}=\left\{\begin{array}{ll}\bar a_{i-2j,j},&\mbox{if }\ i\geq
2j\\0,&\mbox{otherwise}\ea\right\},
\ee
and similarly for $y$.

\section*{Appendix C}

\setcounter{section}{3}\setcounter{equation}{0}

\noindent
First, we conf\/irm that the solution (4.2), namely,
\be
\ba{l}
\ds x=\sum_{k=0}^{\infty}\sum_{l=0}^{\infty}
\sum_{m=0}^{\infty}a_{klm}\tau^{k+\bar\alpha+(2+\bar\alpha)m+\bar rl},
\vspace{3mm}\\
\ds y=\sum_{k=0}^{\infty}\sum_{l=0}^{\infty}\sum_{m=0}^{\infty}b_{klm}
\tau^{k-2+(2+\bar\alpha)m+\bar rl},
\ea
\ee
of the system (1.4), corresponding to the case(ii)-singularities, is
self-consistent and contains four arbitrary constants.

Substitution of (C1) into (1.4) gives
\be
a_{000}\quad \mbox{ is arbitrary},\qquad  b_{000}=6,
\ee
and for $(k,l,m)\neq(0,0,0)$, the coef\/f\/icient recursion relations
\be
\ba{l}
\ds [k+(2+\bar\alpha)m+\bar r l][k+(2+\bar\alpha)m+\bar r(l-1)]a_{klm}+
2\lambda a_{000}b_{klm}\vspace{3mm}\\
\ds \qquad =-Aa_{k-2,l,m}-2\lambda\sum_{(0,0,0)\prec(p,q,r)
\prec(k,l,m)}a_{k-p,l-q,m-r}b_{pqr},
\vspace{3mm}\\
\ds [k+(2+\bar\alpha)m+\bar r l+1][k+(2+\bar\alpha)m+\bar rl-6]b_{klm}
\vspace{3mm}\\
\ds \qquad=-Bb_{k-2,l,m}-\lambda\sum_{p=0}^{k}\sum_{q=0}^{l}
\sum_{r=0}^{m}a_{k-p,l-q,m-r-2}a_{pqr}\vspace{3mm}\\
\ds \qquad +\sum_{(0,0,0)\prec(p,q,r)\prec(k,l,m)}b_{k-p,l-q,m-r}b_{pqr},
\ea
\ee
where we have def\/ined the ``$\prec$'' symbol as in Appendix~B.
Expressed in matrix form, the left-hand sides become
\be
\hspace*{-10mm}
\left[\!\! \begin{array}{cc}
\mbox{\scriptsize $[k+(2+\bar \alpha)m+\bar r l][k+(2+\bar\alpha)m+
\bar r(l-1)]$}
 &2\lambda a_{000}
 \vspace{2mm}\\
 0& \mbox{\scriptsize $[k+(2+\bar \alpha)m+\bar r
l+1][k+(2+\bar\alpha)m+\bar rl-6]$} \ea\!\!\right]\!\!
\left[\ba{c}a_{klm}\\b_{klm}\ea\right],
\ee
and the coef\/f\/icient matrix is singular precisely when
$k+(2+\bar\alpha)m+\bar r l$ has one of the values $-1, 6, \bar r$,
i.e., when
\be
(k,l,m)=(-1,0,0), (6,0,0), (0,1,0).
\ee
These, together with $(k,l,m)=(0,0,0)$, are the resonances of (C1).
Compatibility at $(0,0,0)$ and the arbitrariness of $a_{000}$ has
already been conf\/irmed. The compatibility of $(-1,0,0)$, ref\/lecting
the arbitrariness of $t_0$, is trivial, as usual. It remains to check
the remaining two.

At $(k,l,m)=(0,1,0)$, we f\/ind that
\be
b_{010}=0,\qquad  a_{010}\quad \mbox{is arbitrary}.
\ee
The compatibility condition at $(k,l,m)=(6,0,0)$ is
\be
-Bb_{400}+\sum_{p=1}^5 b_{6-p,0,0}b_{p00}=0.
\ee
The coef\/f\/icients which af\/fect (C7) are
\be
b_{100}=0,\qquad b_{300}=0,\qquad b_{500}=0,\qquad
b_{200}=\frac{B}{2},\qquad b_{400}=\frac{B^2}{40},
\ee
from which (C7) is conf\/irmed. The second equation in (C3) then shows
that $b_{600}$ is arbitrary, and the f\/irst gives $a_{600}$ in terms
of $b_{600}$.

Second, we note that a solution dif\/ferent from (C1) has been given in
[3], but we f\/ind that the coef\/f\/icient recursion relations (their
equation (3.9)) do not follow from the solution (their equation
(3.6)) without neglecting certain terms, precisely those terms which
lead to anomalies when a resummation is performed in an attempt to
prove convergence.

Finally, we remark upon the ``derivation'' of the form of (C1). The
exponent $\bar rl$ is required in order to produce an arbitrary
coef\/f\/icient corresponding to the resonance at $\bar r$. This is
usual. The exponent $\bar \alpha m$ is required in order that all terms
in the system (1.4) can be balanced. But since $Re(\bar\alpha)>-2$,
$\bar\alpha m$ must be replaced by $(2+\bar\alpha)m$ in order to avoid
exponents of $\tau$ with arbitrarily large negative real parts, which
would render the series invalid.

\label{melkonian-lp}
\end{document}